\documentclass[11pt]{amsart}

\usepackage{pb-diagram}
\usepackage{amssymb}

\newtheorem{thm}{Theorem}[section]   
\newtheorem{cor}[thm]{Corollary}     
\newtheorem{lem}[thm]{Lemma}         
\newtheorem{prop}[thm]{Proposition}  

\theoremstyle{definition}

\theoremstyle{remark}
\newtheorem{rem}[thm]{Remark}        

\numberwithin{equation}{section}     

\newcommand{\secref}[1]{Section~\textup{\ref{#1}}}
\newcommand{\thmref}[1]{Theorem~\textup{\ref{#1}}}
\newcommand{\corref}[1]{Corollary~\textup{\ref{#1}}}
\newcommand{\lemref}[1]{Lemma~\textup{\ref{#1}}}
\newcommand{\propref}[1]{Proposition~\textup{\ref{#1}}}

\newcommand{\midtext}[1]{\quad\text{#1}\quad}
\newcommand{\righttext}[1]{\qquad\text{#1 }}

\DeclareMathOperator{\id}{id}

\DeclareMathOperator{\ind}{Ind}
\DeclareMathOperator{\Ind}{Ind}
\DeclareMathOperator{\Aut}{Aut}

\DeclareMathOperator{\supp}{supp}

\DeclareMathOperator*{\clsp}{\overline{span}}
\DeclareMathOperator*{\spn}{span}

\renewcommand{\c}[1]{\mathcal #1}
\renewcommand{\b}[1]{\mathbb #1}

\newcommand{\cross}{\times}

\newcommand{\norm}[1]{\left\| #1 \right\|}
\newcommand{\abs}[1]{\left| #1 \right|}

\newcommand{\case}[2]{\quad\text{if $#1 = #2$ (and $0$ if not)}}
\newcommand{\twocase}[4]{\quad\text{if $#1 = #2$ and $#3 = #4$ 
                         (and $0$ if not)}}
\newcommand{\inner}[1]{\left\langle #1 \right\rangle}

\newcommand{\rip}[2]{
  \inner{#1}_{\!{#2}}  }

\newcommand{\what}{\widehat}
\newcommand{\deltahat}{\what\delta}

\newcommand{\slashstyle}{\mathcal}

\newcommand{\B}{{\slashstyle B}}

\newcommand{\K}{{\slashstyle K}}

\newcommand{\ess}{{s}}
\newcommand{\pee}{{p}}
\newcommand{\tee}{{t}}
\newcommand{\cue}{{q}}

\renewcommand{\forall}{\righttext{for}}

\begin{document}

\title{Skew products and crossed
products by coactions}

\author[Kaliszewski]{S.~Kaliszewski}
\address{Department of Mathematics \\ Arizona State University\\
Tempe, AZ 85287}
\email{kaz@math.la.asu.edu}

\author[Quigg]{John Quigg}
\email{quigg@math.la.asu.edu}

\author[Raeburn]{Iain Raeburn}
\address{Department of Mathematics \\ University of Newcastle \\
NSW 2308 \\ Australia}
\email{iain@frey.newcastle.edu.au}

\thanks{Research partially supported by National Science Foundation
Grant DMS9401253 and the Australian Research Council}




\begin{abstract}
Given a labeling
$c$ of the edges of a directed graph $E$ by elements of a 
discrete group $G$, one can form a skew-product graph $E\cross_c G$. 
We show, using the universal properties of the various constructions
involved, that there is
a coaction $\delta$ of $G$ on $C^*(E)$ such that $C^*(E\cross_c G)$ is
isomorphic to the crossed product $C^*(E)\cross_\delta G$. 
This isomorphism is equivariant for the dual action $\what\delta$ and a
natural action $\gamma$ of $G$ on $C^*(E\times_c G)$; 
following results of Kumjian and Pask, we show that
$$ C^*(E\times_c G)\times_\gamma G \cong C^*(E\times_c
G)\times_{\gamma,r}G \cong C^*(E)\otimes\K(\ell^2(G)),
\phantom{xxxxxxxxxxxxx}$$
and it turns out that the action $\gamma$ is always amenable.
We also obtain 
corresponding results for $r$-discrete groupoids $Q$ and continuous
homomorphisms $c\colon Q\to G$, provided $Q$ is amenable.  
Some of these hold under a more general technical condition
which obtains whenever $Q$ is amenable or second-countable. 
\end{abstract}

\maketitle


\section{Introduction}
\label{intro-sec}

The $C^*$-algebra of a directed graph $E$ is the universal
$C^*$-algebra $C^*(E)$ generated by a family of partial isometries which
are parameterized by the edges of the graph and satisfy relations of
Cuntz-Krieger type reflecting the structure of the graph. A
labeling $c$ of the edges by elements of a discrete group $G$ gives
rise to a skew-product graph $E\times_c G$, and the natural
action of $G$ by translation on $E\times_c G$ lifts to an action
$\gamma$ of $G$ by automorphisms of
$C^*(E\times_c G)$. Kumjian and Pask have
recently proved (\cite[Corollary 3.9]{KP-CD}) that 
\begin{equation}
\label{KP-eq}
C^*(E\times_c G)\times_\gamma G\cong 
C^*(E)\otimes\K(\ell^2(G)).
\end{equation}
{From} this they
obtained an elegant description of the crossed product 
$C^*(F)\times_\beta G$
arising from a free action of $G$ on a graph $F$ (\cite[Corollary
3.10]{KP-CD}).

Kumjian and Pask studied $C^*(E\times_c G)$ by  observing that the
groupoid model for $E\times_c G$ is a skew product of the groupoid
model for $E$, and establishing an analogous stable isomorphism
for the
$C^*$-algebras of skew-product groupoids. They also mentioned that one
could obtain these stable isomorphisms from duality theory and a result of
Masuda (see \cite[Note 3.7]{KP-CD}). This second argument raises some
interesting issues, which are settled in this paper.

We begin in \secref{skew-graph-sec} 
by tackling graph $C^*$-algebras directly. We show that
$C^*(E\times_c G)$ can be realized as the crossed product
$C^*(E)\times_\delta G$ by a coaction $\delta$ of $G$ (see
Theorem~\ref{eqvt-isom}), and apply the duality theorem of Katayama
\cite{Kat-TD} to deduce that
\begin{equation}\label{kpiso}
C^*(E\times_c G)\times_{\gamma,r}G\cong(C^*(E)\times_\delta
G)\times_{\hat\delta,r}G\cong C^*(E)\otimes\K(\ell^2(G))
\end{equation}
(see Corollary~\ref{red-stable-cor}). Since Katayama's
theorem involves the reduced crossed product, the result in
(\ref{kpiso}) is slightly different from Kumjian and Pask's
(\ref{KP-eq})
concerning full crossed products. Together, these two results
suggest that the action of $G$ on
$C^*(E\times_c G)$ should be amenable; we prove this in Section 3 by
giving a new proof of the Kumjian-Pask theorem which allows us to see
directly 
that the regular representation of $C^*(E\times_c G)\times_\gamma G$ is
faithful. 

Our proof of the Kumjian-Pask theorem is elementary in the sense that
it uses only the universal properties of graph $C^*$-algebras, and avoids
groupoid and other models. It is therefore slightly more general,
and will appeal to those who are primarily interested in graph
$C^*$-algebras. Aficionados of groupoids, however, will naturally ask if we
can produce similar results for the
$C^*$-algebra
$C^*(Q\times_c G)$ of a skew-product groupoid $Q\times_c G$. We do
this (at least for $r$-discrete groupoids $Q$) in the second half of
the paper.

Masuda has already identified the groupoid algebra $C^*(Q\times_c G)$ as
a crossed product by a coaction, in the context of spatially-defined
groupoid $C^*$-algebras, coactions and crossed products (\cite[Theorem
3.2]{Mas-GD}). Nowadays, one would prefer to use full coactions and
crossed products, and to give arguments which exploit their universal
properties. The result we obtain this way, \thmref{gpdiso},
is more general than could be
deduced from \cite{Mas-GD}, and highlights an intriguing technical
question: does the $C^*$-algebra of the subgroupoid
$c^{-1}(e)=\{q\in Q \mid c(q)=e\}$ embed faithfully in $C^*(Q)$? We 
answer this in the affirmative 
for $Q$ amenable (\lemref{amen-cond-lem}) or second countable
(\thmref{faith}).

In Section 5, we establish the amenability of the canonical action of
$G$ on $C^*(Q\times_c G)$ when $Q$ is amenable. The results of Section
5 are analogous to those of Section 3, but here we show directly that
the action is amenable (\propref{amenact}) using the theory of
\cite{AR-AG} and \cite{MRW-EI}, and deduce the original version
of
\cite[Theorem 3.7]{KP-CD} for full crossed products. 

\subsection*{Conventions.} A \emph{directed graph} is a quadruple
$E=(E^0,E^1,r,s)$ consisting of a set $E^0$ of vertices, a set $E^1$
of edges and maps $r,s\colon E^1\to E^0$ describing the range and source of
edges. (This notation is becoming standard because one can then write
$E^n$ for the set of paths of length $n$, and think of vertices as
paths of length $0$.) The graph $E$ is \emph{row-finite} if each
vertex emits at most finitely many edges. Our graphs may have
sources and sinks.

All groups in this paper are discrete. A \emph{coaction} of a group
$G$ on a $C^*$-algebra $A$ is an injective
nondegenerate homomorphism $\delta$
of $A$ into the spatial tensor product $A\otimes C^*(G)$ such that
$(\delta\otimes\id)\circ\delta=(\id\otimes\delta_G)\circ \delta$. The
\emph{crossed product} $A\times_\delta G$ is the universal $C^*$-algebra
generated by a covariant representation $(j_A,j_G)$ of
$(A,G,\delta)$. In general, we use the conventions of \cite{QuiDC}.
We shall write $\lambda$ for the left regular representation of a
group $G$ on
$\ell^2(G)$, $\rho$ for the right regular representation, and $M$ for
the representation of $C_0(G)$ by multiplication operators on
$\ell^2(G)$.  The characteristic function of a set $K$ will be denoted
by~$\chi_K$.  


\section{Skew-product graphs and duality}
\label{skew-graph-sec}

Let 
$E=(E^0,E^1,r,s)$ be a row-finite directed graph.  Following
\cite{KPR-CK},  a \emph{Cuntz-Krieger
$E$-family} is a collection $\{ \tee_f, \cue_v\mid f\in E^1, v\in E^0\,\}$
of partial isometries $\tee_f$ and mutually orthogonal projections
$\cue_v$ in a $C^*$-algebra $B$ such that
\begin{equation*}
\tee_f^*\tee_f = \cue_{r(f)} \midtext{and}
\cue_v =\sum_{s(f)=v} \tee_f\tee_f^*
\end{equation*}
for each $f\in E^1$ and every $v\in E^0$ which is not a sink.  
By \cite[Theorem~1.2]{KPR-CK}, there is an essentially unique
$C^*$-algebra
$C^*(E)$, generated
by a Cuntz-Krieger $E$-family $\{\,\ess_f, \pee_v\,\}$,
which is universal in the sense that for any Cuntz-Krieger
$E$-family  $\{\, \tee_f, \cue_v\,\}$ in a $C^*$-algebra $B$, there is a
homomorphism $\Phi=\Phi_{\tee,\cue}\colon C^*(E)\to B$ such that 
$\Phi(\ess_f)=\tee_f$ and $\Phi(\pee_v)=\cue_v$ for 
$f\in E^1, v\in E^0$.
If $\sum_v \cue_v\to 1$ strictly in $M(B)$,  we say that $\{\, \tee_f,
\cue_v\,\}$ is a \emph{nondegenerate} Cuntz-Krieger $E$-family, and the
homomorphism $\Phi_{\tee,\cue}$ is then nondegenerate. 
Products $s_e^*s_f$ cancel (see \cite[Lemma~1.1]{KPR-CK}),
so $C^*(E)$ is densely spanned by 
the projections $\pee_v$ and
products of the form
$\ess_\mu \ess_\nu^*=\ess_{e_1}\ess_{e_2}\dots
\ess_{e_n}\ess_{f_m}^*\dots\ess_{f_1}^*$,
where $\mu$ and $\nu$ are finite paths in the graph $E$.

For each $z\in {\b T}$,
$\{\, z\ess_f, \pee_v\,\}$ is a Cuntz-Krieger $E$-family, so there is an
automorphism $\alpha_z$ of
$C^*(E)$ such that
$\alpha_z(\ess_f) = z\ess_f$ and $\alpha_z(\pee_v) = \pee_v$.
For each pair
of paths $\mu,\nu$ the  map $z\mapsto \alpha_z(\ess_\mu \ess_\nu^*)$ is
continuous, and it follows from a routine $\epsilon/3$-argument
that
$\alpha$ is a strongly continuous action of
$\b T$ on
$C^*(E)$. It was observed in
\cite{aHR-IS} that the existence of this \emph{gauge action} $\alpha$
characterizes the universal $C^*$-algebra 
$C^*(E)$. The following extension of \cite[Theorem~2.3]{aHR-IS} will
appear in \cite{BPR}; it is proved by modifying the proof of
\cite[Theorem~2.3]{aHR-IS} to allow for infinite graphs and the
possibility of sinks.
 
\begin{lem}\label{gauge-lem}
Let $E$ be a row-finite directed graph, and suppose $B$ is a $C^*$-algebra
generated by a Cuntz-Krieger $E$-family $\{\, \tee_f, \cue_v\,\}$.  If
all the $\tee_f$ and $\cue_v$ are non-zero and there is a strongly continuous
action
$\beta$ of
$\b T$ on
$B$ such that
$\beta_z(\tee_f) = z\tee_f$ and $\beta_z(\cue_v)=q_v$, then the canonical
homomorphism $\Phi_{\tee,\cue}\colon C^*(E)\to B$ is an isomorphism.  
\end{lem}

A \emph{labeling} of $E$ by a
group $G$ is just a function $c\colon E^1\to G$.
The \emph{skew-product} graph $E\times_c G$ is the directed
graph with vertex set $E^0\times G$, edge set $E^1\times G$, and range
and source maps defined by
\[
r(f,t) = (r(f),t)\midtext{and} s(f,t) = (s(f),c(f)t)
\forall  (f,t)\in E^1\times G.
\]
Since
$s^{-1}(v,t) = \{\,(f,c(f)t)\mid f\in s^{-1}(v)\,\}$, the vertex
$(v,t)\in(E\times_cG)^0$ emits the same number of edges as $v\in E^0$;
thus
$E\times_cG$ is row-finite if and only if $E$ is, and
$(v,t)$ is a sink if and only if $v$ is.

\begin{rem}
Our skew product $E\times_cG$
is not quite the same as the versions $E(c)$ in  
\cite{KP-CD} and $E^c$ in \cite{GT-TG};
however, there are isomorphisms
$\phi\colon E(c)\to E\times_cG$ 
and $\psi\colon E^c\to E\times_cG$ given by
\[
\phi(t,v) = \psi(v,t) = (v,t^{-1})
\midtext{and}
\phi(t,f) = \psi(f,t) = (f,c(f)^{-1}t^{-1}).
\]
Our conventions were chosen to make the isomorphism of
\thmref{eqvt-isom} more natural.
\end{rem}

\begin{lem}\label{is-delta}
Let $c$ be a labeling of a row-finite directed graph $E$ by a
discrete group $G$.  
Then there is a coaction $\delta$ of $G$ on $C^*(E)$ such that
\begin{equation}\label{delta-eq}
\delta(\ess_f) = \ess_f\otimes c(f)\midtext{and}
\delta(\pee_v) = \pee_v\otimes 1\righttext{for} f\in E^1, 
\mbox{ } v\in E^0.
\end{equation}
\end{lem}

\begin{proof}
Straightforward calculations show that $\{\,
\ess_f\otimes c(f),\, \pee_v\otimes 1\,\}$ is a 
nondegenerate Cuntz-Krieger $E$-family, so the
universal property gives a 
nondegenerate homomorphism $\delta\colon C^*(E)\to
C^*(E)\otimes C^*(G)$ which satisfies \eqref{delta-eq}.
\lemref{gauge-lem} implies that $\delta$ is injective: take
$\beta=\alpha\otimes\id$, where
$\alpha$ is the gauge action of $\b T$ on $C^*(E)$.  It follows from
\eqref{delta-eq} that the coaction identity
$(\delta\otimes\id)\circ\delta = (\id\otimes\delta_G)\circ\delta$ holds
on generators $\ess_f$ and $\pee_v$, and it extends by algebra and continuity
to all of
$C^*(E)$.  
\end{proof}

The group  $G$ acts on the graph
$E\times_cG$ by right translation, so that $t\cdot(v,s)=
(v,st^{-1})$  and $t\cdot(f,s)= (f,st^{-1})$; this induces an
action $\gamma:G\to\Aut C^*(E\times_c G)$ such that
\begin{equation}\label{gamma-eq}
\gamma_t(\ess_{(f,s)}) = \ess_{(f,st^{-1})}\midtext{and}
\gamma_t(\pee_{(v,s)}) = \pee_{(v,st^{-1})}.
\end{equation}

\begin{thm}\label{eqvt-isom}
Let $c$ be a labeling of a row-finite directed graph $E$
by a discrete group $G$, and let $\delta$ 
be the coaction from \lemref{is-delta}. 
Then 
$$C^*(E)\times_\delta G\cong C^*(E\times_cG),$$
equivariantly for the dual action $\what\delta$ and 
the action $\gamma$ of Equation~\eqref{gamma-eq}.
\end{thm}

\begin{proof}
We use the calculus of \cite{EQ-IC} to handle elements of the
crossed product $C^*(E)\times_\delta G$.  
For each $t\in G$, let $C^*(E)_t =
\{\, a\in C^*(E)\mid \delta(a) = a\otimes t\,\}$ denote the 
corresponding spectral
subspace; we write $a_t$ to denote a
generic element of $C^*(E)_t$.
(This subscript convention conflicts with the
standard notation for Cuntz-Krieger families: each partial isometry
$\ess_f$ is in $C^*(E)_{c(f)}$, and each projection $\pee_v$ is in $C^*(E)_e$,
where $e$ is the identity element of $G$. 
We hope this does not cause confusion.)
Then $C^*(E)\times_\delta G$ is densely spanned by the set $\{\,
(a_t,u)\mid a_t\in C^*(E)_t;\, t,u\in G\,\}$, and the algebraic
operations are given on this set by
\[
(a_r,s)(a_t,u) = (a_ra_t,u)\case{s}{tu},
\mbox{\ and\quad} (a_t,u)^* = (a_t^*, tu).
\]
(If $(j_{C^*(E)},j_G)$ is the canonical covariant homomorphism of
$(C^*(E),C_0(G))$ into
$M(C^*(E)\times_\delta G)$, then $(a_t,u)$ is by
definition $j_{C^*(E)}(a_t)j_G(\chi_{\{u\}})$.)  
The dual action $\what\delta$ of
$G$ on $C^*(E)\times_\delta G$ is characterized by
$\what\delta_s(a_t,u) = (a_t,us^{-1})$.

We aim to define a Cuntz-Krieger $E\times_cG$-family
$\{\tee_{(f,t)},\cue_{(v,t)}\}$ in $C^*(E)\times_\delta G$ by putting
\[
\tee_{(f,t)} = (\ess_f,t)\midtext{and}  \cue_{(v,t)} = (\pee_v,t)
\]  
for $(f,t)\in (E\times_cG)^0$ and $(v,t)\in (E\times_cG)^1$.
To see that this is
indeed a Cuntz-Krieger family, note first that $\pee_v\in
C^*(E)_e$ for all vertices $v$, so the 
$\cue_{(v,t)}$ are mutually orthogonal
projections.  
Next note that $\ess_f\in C^*(E)_{c(f)}$, so
\[
\tee_{(f,t)}^*\tee_{(f,t)} 
=(\ess_f^*,c(f)t)(\ess_f,t)= (\ess_f^*\ess_f,t)=(\pee_{r(f)},t)=\cue_{r(f,t)};
\]
if $(v,t)$  is a not sink, then
$v$ is not a sink in $E$, so
\begin{eqnarray*}
 \cue_{(v,t)}
& = &
(\pee_v,t)=\sum_{s(f)=v}(\ess_f\ess_f^*,t)=
\sum_{s(f)=v}(\ess_f,c(f)^{-1}t)(\ess_f^*,t)\\&=&
\sum_{s(f)=v}(\ess_f,c(f)^{-1}t)(\ess_f,c(f)^{-1}t)^*=
\sum_{s(f,r)=(v,t)}\tee_{(f,r)}\tee_{(f,r)}^*.
\end{eqnarray*}
This shows that $\{\tee_{(f,t)},\cue_{(v,t)}\}$ is a Cuntz-Krieger
$E\times_cG$-family. 

The universal property of the graph algebra
now gives a homomorphism
$\Phi=\Phi_{\tee,\cue}\colon C^*(E\times_cG)\to C^*(E)\times_\delta G$ such
that
$\Phi(\ess_{(f,t)}) = \tee_{(f,t)}$ and
$\Phi(\pee_{(v,t)}) =  \cue_{(v,t)}$;
we shall prove that it is an isomorphism using \lemref{gauge-lem}. The
gauge action $\alpha$ of
$\b T$ on
$C^*(E)$ commutes with the coaction $\delta$, in the sense that
$\delta(\alpha_z(a))= \alpha_z\otimes\id(\delta(a))$ for each $z\in {\b
T}$ and $a\in C^*(E)$; it therefore induces an action
$\alpha\times \id$ of $\b T$ on $C^*(E)\times_\delta G$ such that
\[
(\alpha\times \id)_z(\tee_{(f,t)})=(\alpha\times \id)_z(\ess_f,t) = (z\ess_f,t) =
z\tee_{(f,t)}\midtext{and} (\alpha\times \id)_z( \cue_{(v,t)})= q_{(v,t)}.
\]
One can see that the elements of $\{\tee_{(f,t)},\cue_{(v,t)}\}$ are all
non-zero by fixing a faithful representation $\pi$ of $C^*(E)$ and
considering the regular representation
$\Ind\pi=((\pi\otimes\lambda)\circ\delta)\times(1\otimes M)$ of
$C^*(E)\times_\delta G$ induced by $\pi$: the operator
$\Ind\pi(\tee_{(f,t)})$, for example, is just $(\pi(\ess_f)\otimes
\lambda_{c(f)})(1\otimes M(\chi_{\{t\}}))$, which has non-zero initial
projection $\pi(\pee_{r(f)})\otimes M(\chi_{\{t\}})$. Since 
$(\ess_e, c(f)t)(\ess_f, t) = (\ess_e\ess_f, t)$ and
$(\ess_e,c(f)^{-1}t)(\ess_f^*, t) = (\ess_e\ess_f^*, t)$, 
the range of $\Phi$
contains the generating family 
$\{j_{C^*(E)}(\ess_\mu\ess_\nu^*)j_G(\chi_{\{t\}}),
j_{C^*(E)}(\pee_v)j_G(\chi_{\{t\}})\}$,
and hence is all of $C^*(E)\times_\delta G$.
Thus \lemref{gauge-lem} applies, and
$\Phi$ is an isomorphism of
$C^*(E\times_cG)$ onto $C^*(E)\times_\delta G$.  

Finally, we check that $\Phi$ intertwines $\gamma$ and $\what\delta$:
\[
\Phi(\gamma_r(\ess_{(f,t)})) = 
\Phi(\ess_{(f,tr^{-1})})=(\ess_f,tr^{-1})=\what\delta_r(\ess_f,t)=
\what\delta_r(\Phi(\ess_{(f,t)})),
\]
and this completes the proof.
\end{proof}

\begin{cor}\label{red-stable-cor}
Let  $c$ be a labeling of a row-finite directed graph $E$
by a discrete group $G$, and let $\gamma$ be the action 
of Equation~\eqref{gamma-eq}.  Then
$$C^*(E\times_cG)\times_{\gamma,r}G \cong C^*(E)\otimes\K(\ell^2(G)).$$
\end{cor}

\begin{proof}
The corollary follows from \thmref{eqvt-isom} 
and Katayama's duality theorem \cite[Theorem~8]{Kat-TD}. (Even
though we are using full coactions, Katayama's theorem still applies:
the regular representation is an isomorphism
of $C^*(E)\times_\delta G$ onto the (reduced) crossed product by the
reduction of $\delta$; see \cite[Corollary~2.6]{NilDF}).
\end{proof}


\section{Skew-product graphs: the full crossed product}
\label{full-graph-sec}

\begin{thm}\label{direct-isom-thm}
Let $c$ be a labeling of a row-finite directed graph $E$
by a discrete group $G$, and let $\gamma$ be the action of
$G$ defined by Equation~\eqref{gamma-eq}.  Then
$$C^*(E\times_cG)\times_\gamma G\cong C^*(E)\otimes\K(\ell^2(G)).$$
\end{thm}

\begin{proof}
Since $G$ is discrete, 
$C^*(E\times_cG)\times_\gamma G$ is
generated by the set of products 
$\{\, \ess_{(f,r)}u_t,\, \pee_{(v,r)}u_t\,\}$, where 
$\{\,\ess_{(f,r)},\, \pee_{(v,r)}\,\}$ is a 
nondegenerate Cuntz-Krieger
$E\times_cG$-family and $u$ is the canonical homomorphism of
$G$ into $U\!M(C^*(E\times_cG)\times_\gamma G)$ satisfying 
\begin{equation}\label{univ1-eq}
u_t \ess_{(f,r)} = \ess_{(f,rt^{-1})}u_t\midtext{and}
u_t \pee_{(v,r)} = \pee_{(v,rt^{-1})}u_t\righttext{for}t\in G.
\end{equation}
Moreover, the crossed product is universal in the sense that for any
nondegenerate Cuntz-Krieger $E\times_cG$-family $\{\,\tee_{(f,r)},
\cue_{(v,r)}\,\}$ in a $C^*$-algebra $B$ and any homomorphism $v$ of $G$
into $U\!M(B)$ satisfying the analogue of \eqref{univ1-eq}, there is a
unique nondegenerate homomorphism $\Theta =\Theta_{\tee,\cue,v}$ of 
$C^*(E\times_cG)\times_\gamma G$ into $B$
which takes each generator to its counterpart in $B$.  

We now construct such a family 
$\{\, \tee_{(f,r)}, \cue_{(v,r)}, v_t\,\}$ 
in $C^*(E)\otimes \K(\ell^2(G))$. With
$\{\ess_f,\pee_v\}$ denoting the canonical generators of $C^*(E)$ and writing
$\chi_r$ for $M(\chi_{\{r\}})$, we set 
\begin{equation*}
\tee_{(f,r)} = \ess_f\otimes \lambda_{c(f)}\chi_r\midtext{and} 
\cue_{(v,r)} = \pee_v\otimes \chi_r.
\end{equation*}
Then the $\cue_{(v,r)}$ are clearly mutually orthogonal projections, 
and $\sum_{v,r} \cue_{(v,r)}\to 1$ strictly in
$M(C^*(E)\otimes\K(\ell^2(G)))$. Further, we have
\[
\tee_{(f,r)}^*\tee_{(f,r)}
=\ess_f^*\ess_f\otimes
\chi_r^*\lambda_{c(f)}^*\lambda_{c(f)}\chi_r=\ess_f^*\ess_f\otimes
\chi_r=\pee_{r(f)}\otimes \chi_r= \cue_{r(f,r)},
\]
and
\begin{eqnarray*}
\cue_{(v,r)}
& = & \sum_{s(f)=v}\ess_f\ess_f^*\otimes \chi_r\\
& = & \sum_{s(f)=v}(\ess_f\otimes \chi_r\lambda_{c(f)})(\ess_f\otimes
\chi_r\lambda_{c(f)})^*\\
& = & \sum_{s(f)=v}(\ess_f\otimes \lambda_{c(f)}\chi_{c(f)^{-1}r})
(\ess_f\otimes
\lambda_{c(f)}\chi_{c(f)^{-1}r})^*\\
& = & \sum_{s(f)=v}\tee_{(f,c(f)^{-1}r)}\tee_{(f,c(f)^{-1}r)}^*\\
& = & \sum_{s(f,t)=(v,r)} \tee_{(f,t)}\tee_{(f,t)}^*,
\end{eqnarray*}
so $\{\tee_{(f,r)},\cue_{(v,r)}\}$ is a Cuntz-Krieger $E\times_c G$-family.
The unitary elements $1\otimes \rho_t$ of $M(C^*(E)\otimes
\K(\ell^2(G)))$ satisfy
\begin{eqnarray*}
(1\otimes \rho_t)\tee_{(f,r)} = \ess_f\otimes
\rho_t\lambda_{c(f)}\chi_r &=&\ess_f\otimes
\lambda_{c(f)}\chi_{rt^{-1}}\rho_t = \tee_{(f,rt^{-1})}(1\otimes \rho_t), \
\mbox{ and}\\
(1\otimes \rho_t)\cue_{(v,r)} = \pee_v\otimes \rho_t\chi_r &=& \pee_v\otimes
\chi_{rt^{-1}}\rho_t = \cue_{(v,rt^{-1})}(1\otimes \rho_t);
\end{eqnarray*}
thus we get a 
nondegenerate homomorphism $\Theta=\Theta_{\tee,\cue,1\otimes\rho}\colon
C^*(E\times_cG)\times_\gamma G\to C^*(E)\otimes \K(\ell^2(G))$ such that
\begin{equation}\label{defTheta}
\Theta(\ess_{(f,r)}) = \tee_{(f,r)},\ \ 
\Theta(\pee_{(v,r)}) = \cue_{(v,r)},\ \mbox{ and }
\ \Theta(u_t) = 1\otimes
\rho_t.
\end{equation}  

To construct the inverse for $\Theta$, we use a universal property of
$C^*(E)\otimes \K(\ell^2(G))$. 
Let $\sigma$ denote the action of $G$ on
$C_0(C)$ by right translation: $\sigma_s(f)(t)=f(ts)$. The regular
representation
$M\times
\rho$ is an isomorphism of the crossed product $C_0(G)\times_\sigma
G$ onto $\K(\ell^2(G))$, so we can view $\K(\ell^2(G))$ as the
universal $C^*$-algebra generated by 
the set of products $\{\chi_r\rho_t\mid r,t\in G\}$, where 
$\rho$ is a unitary homomorphism $\rho$ of $G$ and
$\{\,\chi_r\,\}$ is a set 
of mutually orthogonal projections satisfying
\begin{equation}\label{yu-eq}
\rho_t\chi_r = \chi_{rt^{-1}}\rho_t.
\end{equation}
Thus
to get a homomorphism defined on $C^*(E)\otimes \K(\ell^2(G))$ we need
a Cuntz-Krieger $E$-family $\{\tee_f,\cue_v\}$ and a family $\{y_r,u_t\}$ 
analogous to $\{\chi_r, \rho_t\}$ which 
commutes with the Cuntz-Krieger family. 

We begin by constructing a family $\{y_r,u_t\}$ in
$M(C^*(E\times_cG)\times_\gamma G)$. We claim that, for fixed $r\in G$,
the sum $\sum_v \pee_{(v,r)}$ converges strictly in
$M(C^*(E\times_cG)\times_\gamma G)$. Because the canonical
embedding $j_{C^*(E\times_cG)}$ has a strictly continuous
extension, it is enough to check that the sum converges strictly
in $M(C^*(E\times_c G))$. Because all the finite sums are projections,
they have norm uniformly bounded by
$1$, and it is enough to check that 
$\left(\sum_v \pee_{(v,r)}\right)\ess_\mu \ess_\nu^*$ and
$\ess_\mu \ess_\nu^*\left(\sum_v \pee_{(v,r)}\right)$ converge for each pair of
paths $\mu,\nu$ in $E\times_c G$; and that 
$\left(\sum_v \pee_{(v,r)}\right)\pee_{(u,t)}$ and
$\pee_{(u,t)}\left(\sum_v \pee_{(v,r)}\right)$
converge for each vertex $(u,t)$ in $E\times_c G$.  
But in each case these sums reduce to
a single term, so this is trivially true.
Thus we may put $y_r = \sum_v \pee_{(v,r)}\in
M(C^*(E\times_cG)\times_\gamma G)$.  

Now $\{\,y_r\mid r\in G\,\}$ is a mutually orthogonal family of
projections, and $\sum_{v,r}\pee_{(v,r)}\to 1$ strictly in
$M(C^*(E\times_cG))$, so $\sum_s y_s\to 1$ strictly in
$M(C^*(E\times_cG)\times_\gamma G)$.  Moreover, 
if $u$ is the canonical homomorphism of $G$ into
$M(C^*(E\times_cG)\times_\gamma G)$, then 
$$u_ty_r = u_t\sum_v \pee_{(v,r)} = \sum_v
\pee_{(v,rt^{-1})} u_t = y_{rt^{-1}}u_t;$$
thus the family $\{\,y_r,u_t\,\}$ satisfies the analogue of (\ref{yu-eq}),
and therefore gives a
nondegenerate homomorphism $y\times u\colon \K(\ell^2(G))\to
M(C^*(E\times_cG)\times_\gamma G)$. This homomorphism extends to
$\B(\ell^2(G))=M(\K(\ell^2(G)))$, and we can define unitaries
$w_t=y\times u(\lambda_t)$ which satisfy  
$w_ty_r = y_{tr}w_t$ and $w_tu_r = u_rw_t$ for each $r,t\in G$.

Arguing as for the $y_r$ shows that, 
for each fixed $v$ and $f$, the sums
$\sum_r\pee_{(v,r)}$ and
$\sum_r \ess_{(f,r)}$ converge strictly
in $M(C^*(E\times_cG))$. 
Thus we may define $\tee_f$ and $\cue_v$ in 
$M(C^*(E\times_cG)\times_\gamma G)$ by
\[
\tee_f = \left(\sum_r \ess_{(f,r)}\right)w_{c(f)^{-1}}
\midtext{and}
\cue_v = \sum_r\pee_{(v,r)}.
\]
Now $\{\cue_v\}$ is a
family of mutually orthogonal projections;
to check the Cuntz-Krieger relations for
$\{\, \tee_f,\,\cue_v\,\}$, first note that
\[
\biggl(\sum_r\ess_{(f,r)}\biggr)^* \biggl(\sum_t\ess_{(f,t)}\biggr)
=\sum_{r,t}\ess_{(f,r)}^*\ess_{(f,t)}
=\sum_r \ess_{(f,r)}^*\ess_{(f,r)}
=\sum_r \pee_{r(f,r)}
=\cue_{r(f)},
\]
so that
\begin{equation}\label{twq-eq}
\tee_f^*\tee_f =
w_{c(f)}\biggl(\sum_r\ess_{(f,r)}\biggr)^*
\biggl(\sum_t\ess_{(f,t)}\biggr)w_{c(f)}^* = w_{c(f)}\cue_{r(f)}w_{c(f)}^*.
\end{equation}
Easy calculations show that
$y_t\cue_v = q_vy_t$ and $u_tq_v=q_vu_t$, 
so each $\cue_v$ commutes with everything in the range of $y\times u$ in 
$M(C^*(E\times_cG)\times_\gamma G)$, and in particular with each
$w_t$; 
thus Equation (\ref{twq-eq}) implies that $t^*_f\tee_f=\cue_{r(f)}$. 
We also have
\begin{eqnarray*}
\cue_v
& = & \sum_r \sum_{s(f,t)=(v,r)} \ess_{(f,t)}\ess_{(f,t)}^*\\
& = & \sum_r \sum_{s(f)=v} \ess_{(f,c(f)^{-1}r)}\ess_{(f,c(f)^{-1}r)}^*\\
& = & \sum_{s(f)=v} \sum_r \ess_{(f,r)}\ess_{(f,r)}^*\\
& = & \sum_{s(f)=v} \biggl(\sum_r \ess_{(f,r)}w_{c(f)^{-1}}\biggr)
\biggl(\sum_t \ess_{(f,t)}w_{c(f)^{-1}}\biggr)^*\\
& = & \sum_{s(f)=t} \tee_f\tee_f^*,
\end{eqnarray*}
and $\sum_v\cue_v = \sum_{v,r}\pee_{(v,r)}\to 1$ strictly in
$M(C^*(E\times_cG)\times_\gamma G)$, so the set $\{\, \tee_f,\cue_v\,\}$ is 
a nondegenerate Cuntz-Krieger $E$-family. 

We have already observed that each $\cue_v$ commutes with  the
range of $y\times u$.  
Further calculations show that
\begin{eqnarray}
y_s\tee_f 
& = & \biggl(\sum_v \pee_{(v,s)}\biggr) \biggl(\sum_r \ess_{(f,r)}\biggr)
w_{c(f)^{-1}}\notag\\
& = & \ess_{(f,c(f)^{-1}s)}w_{c(f)^{-1}}\label{yt=ty}\\
& = & \biggl(\sum_r \ess_{(f,r)}\biggr)\biggl(\sum_v
\pee_{(v,c(f)^{-1}s)}\biggr)w_{c(f)^{-1}}\notag\\
& = & \biggl(\sum_r\ess_{(f,r)}\biggr)y_{c(f)^{-1}s}w_{c(f)^{-1}}\notag\\
& = & \biggl(\sum_r\ess_{(f,r)}\biggr)w_{c(f)^{-1}}y_s\notag\\
& = & \tee_fy_s\notag
\end{eqnarray}
and
\[
u_t\tee_f 
=u_t\biggl(\sum_r \ess_{(f,r)}\biggr)w_{c(f)^{-1}}
=\biggl(\sum_r
\ess_{(f,rs^{-1})}\biggr)u_tw_{c(f)^{-1}}
=\biggl(\sum_r
\ess_{(f,r)}\biggr)w_{c(f)^{-1}}u_t=\tee_fu_t.
\]
Thus the homomorphisms $\Phi_{\tee,\cue}$ of $C^*(E)$ and $y\times
u$ of $\K(\ell^2(G))$ into $M(C^*(E\times_cG)\times_\gamma G)$ have
commuting ranges, and combine to give a homomorphism $\Upsilon$ of
$C^*(E)\otimes\K(\ell^2(G))$ into $M(C^*(E\times_cG)\times_\gamma G)$
such that $\Upsilon(\ess_f\otimes 1) = \tee_f$, $\Upsilon(\pee_v\otimes 1) = \cue_v$,
$\Upsilon(1\otimes \chi_r) = y_r$, and $\Upsilon(1\otimes \rho_t) = u_t$.  
{From} (\ref{yt=ty}) we deduce that
\begin{equation}\label{into-calc1}
\Upsilon(\ess_f\otimes \chi_r\rho_t)
=  \tee_fy_ru_t =\ess_{(f,c(f)^{-1}r)}w_{c(f)^{-1}}u_t
=  \ess_{(f,c(f)^{-1}r)}u_tw_{c(f)^{-1}};
\end{equation}
since this and $\Upsilon(\pee_v\otimes \chi_r\rho_t)
=  \cue_vy_ru_t=  \pee_{(v,r)}u_t$ belong to $C^*(E\times_cG)\times_\gamma
G$, it follows that
$\Upsilon$ maps $C^*(E)\otimes \K(\ell^2(G))$ into
$C^*(E\times_cG)\times_\gamma G$.  

We shall show that $\Theta$ and $\Upsilon$ are inverses of
one another by checking that $\Upsilon\circ\Theta$ is the identity on the
generating set $\{\ess_{(f,r)}u_t,\, \pee_{(v,r)}u_t\}$ for $C^*(E\times_c
G)\times_\gamma G$, and that
$\Theta\circ \Upsilon$ is the identity on a generating set for
$C^*(E)\otimes \K(\ell^2(G))$.  First we note 
that $T\mapsto \Upsilon(1\otimes T)$ is just $y\times u$
on products $\chi_r\rho_t\in \K(\ell^2(G))$, 
so $\Upsilon(1\otimes \lambda_t)=w_t$ by definition of $w_t$. And since
$\Theta$ extends to a strictly continuous map on 
$M(C^*(E\times_cG)\times_\gamma G)$, we have
$$\Theta(y_ru_t) = \Theta\biggl(\sum_v\pee_{(v,r)}u_t\biggr) =
\sum_v\pee_v\otimes \chi_r\rho_t = 1\otimes \chi_r\rho_t,$$
which implies that $\Theta(w_t)
= 1\otimes \lambda_t$ for $t\in G$. 

We can now compute:
\begin{eqnarray*}
\Upsilon\circ\Theta(\ess_{(f,s)}u_t)
& = & \Upsilon(\ess_f\otimes \lambda_{c(f)}\chi_s\rho_t)\\
& = & \biggl(\sum_r\ess_{(f,r)}\biggr)w_{c(f)^{-1}}w_{c(f)}
\biggl( \sum_v\pee_{(v,s)}\biggr)u_t\\
\end{eqnarray*}
\begin{eqnarray*}
\phantom{\Upsilon\circ\Theta(\ess_{(f,s)}u_t)}
& = & \biggl(\sum_r\ess_{(f,r)}\biggr)\biggl( \sum_v\pee_{(v,s)}\biggr)u_t\\
& = & \ess_{(f,s)}u_t
\end{eqnarray*}
and
\[
\Upsilon\circ\Theta(\pee_{(v,s)}u_t)
=\Upsilon(\pee_v\otimes \chi_s\rho_t)
=\biggl(\sum_r\pee_{(v,r)}\biggr)\biggl(\sum_w \pee_{(w,s)}\biggr)u_t
=\pee_{(v,s)}u_t,
\]
which shows that $\Upsilon\circ\Theta$ is the identity.
Using \eqref{into-calc1} gives
\begin{eqnarray*}
\Theta\circ\Upsilon(\ess_f\otimes \chi_r\rho_t) 
& = & \Theta(\ess_{(f,c(f)^{-1}r)}u_tw_{c(f)^{-1}})\\
& = & \ess_f\otimes \lambda_{c(f)}\chi_{c(f)^{-1}r}\lambda_{c(f)^{-1}}\rho_t\\
& = & \ess_f\otimes \chi_r\rho_t
\end{eqnarray*}
and
\begin{equation*}
\Theta\circ\Upsilon(\pee_v\otimes \chi_r\rho_t) = \Theta(\pee_{(v,r)}u_t) =
\pee_v\otimes \chi_r\rho_t,
\end{equation*}
which shows that $\Theta\circ\Upsilon$ is the identity.
\end{proof}

\thmref{direct-isom-thm} and
\corref{red-stable-cor} imply that $C^*(E\times_cG)\times_\gamma G$
and $C^*(E\times_cG)\times_{\gamma,r}G$
are isomorphic $C^*$-algebras,
so it is natural to ask if the action $\gamma$ is amenable in the sense
that the regular representation of the crossed product is faithful.
To see that it is, consider the following diagram:
\begin{equation}\label{reg-diag}
\begin{diagram}
\node{C^*(E\times_cG)\times_\gamma G}
\arrow[2]{e,t}{\rm\thmref{eqvt-isom}}
\arrow{se,t}{\rm\thmref{direct-isom-thm}}
\arrow[2]{s,l}{\begin{smallmatrix}{\rm regular}\\
{\rm representation}\end{smallmatrix}}
\node[2]{C^*(E)\times_\delta G\times_{\deltahat}G}
\arrow[2]{s,r}{\begin{smallmatrix}{\rm regular}\\
{\rm representation}\end{smallmatrix}}\\
\node[2]{C^*(E)\otimes\K}\\
\node{C^*(E\times_cG)\times_{\gamma,r}G}
\arrow[2]{e,b}{\rm\thmref{eqvt-isom}}
\node[2]{C^*(E)\times_\delta G\times_{\deltahat,r}G.}
\arrow{nw,t}{{\rm Katayama}}
\end{diagram}
\end{equation}
Let 
$(j_{C^*(E)},j_G)\colon (C^*(E),C_0(G))\to M(C^*(E)\times_\delta G)$ and 
$(i_{C^*(E)\times G},i_G)\colon (C^*(E)\times_\delta G,G)\to
M(C^*(E)\times_\delta G\times_{\deltahat}G)$ be the canonical maps. 
Inspection of the formulas on page~768 of
\cite{LPRS-RC} shows that composing the regular representation of
$C^*(E)\times_\delta G\times_{\deltahat}G$ with Katayama's
isomorphism (\cite[Theorem~8]{Kat-TD}) gives
\begin{gather*}
i_{C^*(E)\times G}(j_{C^*(E)}(a)) \mapsto \id\otimes\lambda(\delta(a)),
\mbox{  }i_{C^*(E)\times G}(j_G(g)) \mapsto 1\otimes M(g),
\mbox{  } i_G(t)\mapsto 1\otimes\rho_t
\end{gather*}
for $a\in C^*(E)$, $g\in C_c(G)$, and $t\in G$.  
Thus chasing generators in $C^*(E\times_cG)\times_\gamma G$ round the
outside of the upper right-hand triangle in Diagram~\eqref{reg-diag}
yields
\begin{gather*}
\ess_{(f,r)}\mapsto 
i_{C^*(E)\times G}(j_{C^*(E)}(\ess_f)j_G(\chi_r)) \mapsto 
\id\otimes\lambda(\delta(\ess_f))1\otimes \chi_r =  \tee_{(f,r)},\\
\pee_{(v,r)}\mapsto 
i_{C^*(E)\times G}(j_{C^*(E)}(\pee_v)j_G(\chi_r))\mapsto
\id\otimes\lambda(\delta(\pee_v))1\otimes \chi_r = \cue_{(v,r)},\\
u_t\mapsto i_G(t)\mapsto 1\otimes \rho_t.
\end{gather*}
Since this is exactly what the  isomorphism  $\Theta$ from
\thmref{direct-isom-thm} does (see Equation (\ref{defTheta})), the
upper right-hand corner of Diagram~\eqref{reg-diag} commutes.  But the
outside rectangle commutes by general nonsense, so the lower left-hand
corner commutes too.  This proves:

\begin{cor}\label{amen-cor}
Let  $c$ be a labeling of a row-finite directed graph $E$
by a discrete group $G$.
Then the action $\gamma$ of $G$
from Equation~\eqref{gamma-eq} is amenable in the sense
that the regular
representation of $C^*(E\times_cG)\times_\gamma G$
is faithful.
\end{cor}

\begin{cor}\label{reducedKP2}
Let $G$ be a discrete group acting freely on a row-finite 
directed graph $F$, and
let $\beta$ be the action of $G$ on $C^*(F)$ determined by
$\beta_t(\ess_f)=\ess_{t\cdot f}$ and $\beta_t(\pee_v)=\pee_{t\cdot v}$. Then 
the regular representation of
$C^*(F)\times_\beta G$ is faithful, and 
$$C^*(F)\times_\beta G\cong 
C^*(F)\times_{\beta,r}G\cong 
C^*(F/G)\otimes\K(\ell^2(G)).$$
\end{cor}

\begin{proof}
Since $G$ acts freely, there is a
labeling $c\colon (F/G)^1\to G$ and an isomorphism of $F$ onto
$(F/G)\times_cG$ which carries the given action to the action
of $G$ by right translation (\cite[Theorem~2.2.2]{GT-TG}). Thus this
corollary follows by applying Corollaries~\ref{red-stable-cor}
and~\ref{amen-cor} to $E=F/G$. 
\end{proof}


\section{Skew-product groupoids and duality}
\label{gpd}

We will now give groupoid versions of the results in
\secref{skew-graph-sec}. 
Throughout, we consider a
discrete group $G$, and a groupoid $Q$ which is 
$r$-discrete in the modern sense that the range map $r$
is a local homeomorphism (so that counting measures on the sets
$Q^u=r^{-1}(u)$ for $u$ in the unit space $Q^0$ give a Haar system on
$Q$).
In several of the following arguments, we use the fact
that the $C^*$-algebra of  
an $r$-discrete groupoid $Q$ is the enveloping $C^*$-algebra of
$C_c(Q)$; this follows from \cite[Theorems~7.1 and~8.1]{QS-CA}.

Let $c\colon Q\to G$ be a continuous homomorphism.  
The \emph{skew-product groupoid} $Q\times_c G$ is the
set $Q\times G$ with the product topology and operations given for
$(x,y)\in Q^2$ and $s\in G$ by
\[
(x,c(y)s)(y,s)=(xy,s)\midtext{and}
(x,s)^{-1}=(x^{-1},c(x)s).
\]
Since the range map on the skew-product groupoid is thus
given by $r(x,s) = (r(x),c(x)s)$,
$Q\times_cG$ is $r$-discrete whenever $Q$ is.  
The formula
\begin{equation}
\label{gpdact}
s\cdot (x,t)=(x,ts^{-1})
\righttext{for}s,t\in G,\mbox{ }x\in Q
\end{equation}
defines an action of $G$ by automorphisms of the topological groupoid
$Q\times_c G$.
We let $\beta$ denote
the induced action on $C^*(Q\times_c G)$, which satisfies
\begin{equation}
\label{beta}
\beta_s(f)(x,t)=f\bigl(s^{-1}\cdot (x,t)\bigr)=f(x,ts)
\righttext{for}s,t\in G,\mbox{ }f\in C_c(Q\times_c G),\mbox{ }x\in Q.
\end{equation}

\begin{rem}
It is easily checked that the map
$(x,s)\mapsto(x,c(x)^{-1}s^{-1})$
gives a topological isomorphism of Renault's skew product
\cite[Definition~I.1.6]{RenGA}
onto ours,
which transports Renault's action (also used in
\cite[Proposition~3.7]{KP-CD}) into $\beta$.  
Our conventions were chosen to make the isomorphism of \thmref{gpdiso}
more natural.
\end{rem}

For $s\in G$ define
\begin{equation}
\label{bundle}
C_s=\{f\in C_c(Q)\mid\supp f\subseteq c^{-1}(s)\},
\end{equation}
and put $\c C=\bigcup_{s\in G}C_s$. Then with the operations
from $C_c(Q)$, $\c C$ becomes a
$^*$-algebraic bundle (with incomplete fibers) over $G$ in the sense that
$C_sC_t\subseteq C_{st}$ and $C_s^*=C_{s^{-1}}$.
Since $Q$ is the disjoint union of the open sets $\{c^{-1}(s)\}_{s\in G}$,
we have $\spn_{s\in G}C_s = C_c(Q)$,
which we identify with the space
$\Gamma_c(\c C)$ of finitely supported sections of $\c C$.

\begin{lem}
\label{gpdcoact}
Let $c$ be a continuous homomorphism of an $r$-discrete Hausdorff
groupoid $Q$ into a discrete group $G$.  Then 
there is a coaction $\delta$ of $G$ on $C^*(Q)$ such that
\[
\delta(f_s)=f_s\otimes s\righttext{for}s\in G,\mbox{ }f_s\in C_s.
\]
\end{lem}

\begin{proof}
The above formula extends uniquely to a $^*$-homomorphism of $C_c(Q)$
into $C^*(Q)\otimes C^*(G)$. Since $C^*(Q)$ is the enveloping
$C^*$-algebra of $C_c(Q)$, $\delta$
further extends uniquely to a homomorphism of $C^*(Q)$ into
$C^*(Q)\otimes C^*(G)$. The coaction identity obviously holds on the
generators (that is, the elements of the bundle $\c C$), hence on all of
$C^*(Q)$. The homomorphism $\delta$ is nondegenerate, that is,
\[
\clsp\{\delta(C^*(Q))(C^*(Q)\otimes C^*(G))\}
=C^*(Q)\otimes C^*(G),
\]
since $\delta(f_s)(1\otimes s^{-1}t)=f_s\otimes t$.
To see that $\delta$ is injective, let $1_G$ denote the trivial
one-dimensional representation of $G$, and check on the generators
that $(\id\otimes 1_G)\circ\delta=\id$.
\end{proof}

Let $N=c^{-1}(e)$ be the kernel of the homomorphism $c$, which is an
open subgroupoid of $Q$. Since the restriction of a Haar system to an
open subgroupoid gives a Haar system, counting measures give a Haar
system on $N$, so $N$ is an $r$-discrete groupoid.
The inclusion of $C_c(N)$ in $C_c(Q)$ extends to the enveloping
$C^*$-algebras to give a natural
homomorphism $i$ of $C^*(N)$ into $C^*(Q)$.
For our next results we will need to require that $i$ be faithful.  
We have been unable to show that this holds in general, although 
it does hold when $Q$ is amenable (\lemref{amen-cond-lem}), 
and when $Q$ is second countable (\thmref{faith}).  

\begin{thm}
\label{gpdiso}
Let $c$ be a continuous homomorphism of an $r$-discrete Hausdorff
groupoid $Q$ into a discrete group $G$, let $N=c^{-1}(e)$, 
and let $\delta$ be the
coaction from \lemref{gpdcoact}.  Assume that the natural map $i\colon
C^*(N)\to C^*(Q)$ is faithful.  Then
\[
C^*(Q)\times_\delta G\cong C^*(Q\times_c G),
\]
equivariantly for the dual action $\what\delta$ and the action
$\beta$ of Equation~\eqref{beta}.
\end{thm}

\begin{proof}
Let $\c C$ be the $^*$-algebraic bundle over $G$ defined by
Equation \eqref{bundle},
let $\c C\times G$ be the product bundle over $G\times G$ whose fiber
over $(s,t)$ is $C_s\times\{t\}$, and give $\c C\times G$ the algebraic
operations
\[
(f_s,tu)(g_t,u)=(f_sg_t,u)\midtext{and}
(f_s,t)^*=(f_s^*,st).
\]
Then the space $\Gamma_c(\c C\times G)$ of finitely supported
sections becomes a $^*$-algebra, which can be identified with a dense
$^*$-subalgebra of the crossed product $C^*(Q)\times_\delta G$;
the dual action is characterized by $\what\delta_s(f,t) = (f,ts^{-1})$,
for $s,t\in G$ and $f\in \c C$.  

We claim that $C^*(Q)\times_\delta G$ is the enveloping
$C^*$-algebra of $\Gamma_c(\c C\times G)$. 
Since $C^*(Q)$ is the enveloping
$C^*$-algebra of $\Gamma_c(\c C) = C_c(Q)$, 
by \cite[Theorem~3.3]{EQ-IC} it suffices
to show that the unit 
fiber algebra 
$C^*(Q)_e = \{\, f\in C^*(Q)\mid \delta(f) = f\otimes e\,\}$
of the Fell bundle associated to 
$\delta$ is the enveloping $C^*$-algebra of $C_e$.
To see this, first note that $C^*(Q)_e$
is the closure of $C_e$ in $C^*(Q)$, which in turn is just $i(C^*(N))$
because $i$ maps $C_c(N)$ onto $C_e$.  But $C^*(N)$ is the
enveloping $C^*$-algebra of $C_c(N)$.
Since $i$ is assumed to be faithful, it follows that 
$C^*(Q)_e = i(C^*(N))$ is the
enveloping $C^*$-algebra of $C_e = i(C_c(N))$, and this proves the claim. 

Now for each $s,t\in G$ put
$D_{s,t}
=\{f\in C_c(Q\times_c G)\mid\supp f\subseteq c^{-1}(s)\times\{t\}\}$,
so $C_c(Q\times_c G)=\spn_{s,t\in G}D_{s,t}$.
For $f\in \c C$ and $t\in G$ define $\Psi(f,t)\in C_c(Q\times_c
G)$ by
\begin{equation}
\Psi(f,t)(x,u)=f(x)
\qquad\case{t}{u}.
\end{equation}
Then $\Psi$ extends uniquely to a linear bijection $\Psi$ of
$\Gamma_c(\c C\times G)$ onto $C_c(Q\times_c G)$, since it gives
a linear bijection of each fiber $C_s\times\{t\}$ onto the corresponding
fiber $D_{s,t}$. In fact, $\Psi$ is a homomorphism of
$^*$-algebras. 
It is enough to show $\Psi$ preserves
multiplication and involution. For $s,t,u,v,z\in G$, $f_s\in C_s$,
$g_u\in C_u$, and $x\in Q$,
\begin{align*}
&\bigl(\Psi(f_s,t)\Psi(g_u,v)\bigr)(x,z)
\\&\quad=\sum_{r(y,w)=r(x,z)}\Psi(f_s,t)(y,w)
\Psi(g_u,v)\bigl((y,w)^{-1}(x,z)\bigr)
\\&\quad=\sum_{\begin{smallmatrix}{r(y)=r(x)}\\
{c(y)w=c(x)z}\end{smallmatrix}}f_s(y)
\Psi(g_u,v)\bigl((y^{-1},c(y)w)(x,z)\bigr)
&&\case{t}{w}
\\&\quad=\sum_{r(y)=r(x)}f_s(y)\Psi(g_u,v)(y^{-1}x,z)
&&\case{t}{c(y^{-1}x)z}
\\&\quad=\sum_{r(y)=r(x)}f_s(y)g_u(y^{-1}x)
&&\twocase{t}{uz}{v}{z}
\\&\quad=(f_sg_u)(x)
&&\twocase{t}{uv}{v}{z}
\\&\quad=\Psi(f_sg_u,v)(x,z)
&&\case{t}{uv}
\\&\quad=\Psi\bigl((f_s,t)(g_u,v)\bigr)(x,z),
\end{align*}
and for $s,t,u\in G$, $f_s\in C_s$, and $x\in Q$,
\begin{align*}
\Psi(f_s,t)^*(x,u)
&=\overline{\Psi(f_s,t)\bigl((x,u)^{-1}\bigr)}
\\&=\overline{\Psi(f_s,t)(x^{-1},c(x)u)}
\\&=\overline{f_s(x^{-1})}
&&\case{t}{c(x)u}
\\&=f_s^*(x)
&&\case{st}{u}
\\&=\Psi(f_s^*,st)(x,u)
\\&=\Psi\bigl((f_s,t^*\bigr)(x,u).
\end{align*}
It follows that $\Psi$ extends to an isomorphism of
$C^*(Q)\times_\delta G$ onto $C^*(Q\times_cG)$, since these are
enveloping $C^*$-algebras.  

A straightforward calculation shows that $\Psi$ intertwines the actions
$\what\delta$ and $\beta$.
\end{proof}

\begin{rem}
For $Q$ amenable, the isomorphism of 
\thmref{gpdiso} can be deduced from
\cite[Theorem~3.2]{Mas-GD}, although
Masuda does everything spatially, with reduced
coactions, reduced groupoid $C^*$-algebras, and
crossed products represented on Hilbert space.
To see this, note that 
the amenability of the skew product $Q\times_c G$ 
follows from  that of $Q$ by 
\cite[Proposition~II.3.8]{RenGA}, 
and that 
$C^*(Q)\times_\delta G$ is isomorphic to 
the spatial crossed product by the
reduction of $\delta$ according to results in \cite{QuiFR} and \cite{RaeOCP}.  
\end{rem}

\begin{cor}\label{duality}
With the same hypotheses as \thmref{gpdiso},
\[
C^*(Q\times_c G)\times_{\beta,r}G\cong C^*(Q)\otimes\c K(\ell^2(G)).
\]
\end{cor}

\begin{proof}
This follows immediately from \thmref{gpdiso} and Katayama's
duality theorem \cite[Theorem~8]{Kat-TD}.
(See also the parenthetical remark in the proof of
\corref{red-stable-cor}.)
\end{proof}


\section{Skew-product groupoids: the full crossed product}
\label{gpd-full-sec}

In this section we prove a version of \corref{duality} for full crossed
products, from which we can recover 
Proposition~3.7 of \cite{KP-CD}. 
For this, we shall want to relate semidirect-product groupoids to
crossed products.  In general, if a discrete
group $G$ acts on a topological groupoid $R$, the \emph{semidirect-product
groupoid} $R \rtimes G$
is the product space $R\times G$ with the structure
\[
(x,s)(y,t)=(x(s\cdot y),st)
\midtext{and}
(x,s)^{-1}=(s^{-1}\cdot x^{-1},s^{-1})
\]
whenever this makes sense. (This is readily seen to coincide with
Renault's version in \cite[Definition~I.1.7]{RenGA}.) If $R$ is
$r$-discrete and Hausdorff then so is $R \rtimes G$. 

The following result is presumably folklore, but it never hurts to
record groupoid facts.

\begin{prop}
\label{semi cross}
Let $G$ be a discrete group acting on an $r$-discrete
Hausdorff groupoid $R$, and let $\beta$ denote the associated action
on $C^*(R)$. Then
\[
C^*(R \rtimes G)
\cong
C^*(R)\times_\beta G.
\]
\end{prop}

\begin{proof}
For $f\in C_c(R \rtimes G)$ define
$\Phi(f)\in C_c(G,C_c(R)) \subseteq C_c(G,C^*(R))$
by
\[
\Phi(f)(s)(x)=f(x,s)
\righttext{for}s\in G,x\in R.
\]
Then $\Phi$ is a $^*$-homomorphism, since for $f,g\in C_c(R \rtimes
G)$ we have
\begin{align*}
\bigl( \Phi(f)\Phi(g) \bigr)(s)(x)
&=\sum_t \bigl( \Phi(f)(t)
\beta_t(\Phi(g)(t^{-1}s)) \bigr)(x)
\\&=\sum_t \sum_{r(y)=r(x)} \Phi(f)(t)(y)
\beta_t(\Phi(g)(t^{-1}s))(y^{-1}x)
\\&=\sum_{r(y)=r(x)} \sum_t f(y,t)
g(t^{-1}\cdot (y^{-1}x),t^{-1}s)
\\&=\sum_{r(y,t)=r(x,s)} f(y,t) g((y,t)^{-1}(x,s))
\\&=(fg)(x,s) =\Phi(fg)(s)(x)
\end{align*}
and
\begin{gather*}
\Phi(f)^*(s)(x)
=\beta_s(\Phi(f)(s^{-1})^*)(x)
=\Phi(f)(s^{-1})^*(s^{-1}\cdot x)
=\overline{\Phi(f)(s^{-1})(s^{-1}\cdot x^{-1})}\\
=\overline{f(s^{-1}\cdot x^{-1},s^{-1})}
=\overline{f((x,s)^{-1})}
=f^*(x,s)
=\Phi(f^*)(s)(x).
\end{gather*}
Since $C^*(R \rtimes G)$ is the enveloping $C^*$-algebra of
$C_c(R \rtimes G)$, $\Phi$
extends uniquely to a homomorphism $\Phi$ of $C^*(R \rtimes G)$
into $C^*(R)\times_\beta G$.

To show $\Phi$ is an isomorphism, it suffices to find an inverse for
the map $\Phi\colon C_c(R \rtimes G)\to C_c(G,C_c(R))$,
since $C^*(R)\times_\beta G$ is the
enveloping $C^*$-algebra of the $^*$-algebra $C_c(G,C_c(R))$ 
(see, for example, \cite[Lemma~3.3]{EQ-IC}). Given $f\in
C_c(G,C_c(R))$ define $\Psi(f)\in C_c(R \rtimes G)$ by
\[
\Psi(f)(x,s)=f(s)(x).
\]
Since the support of $\Phi(f)$ in $R\times G$ is just the finite union
of compact sets $\{s\}\times\supp f(s)$ as $s$ runs through $\supp f$, 
$\Psi(f)$ has compact support. Moreover, it is obvious that $\Psi$ is 
the required inverse for $\Phi$ at the level of $C_c$-functions. 
\end{proof}

To show that the isomorphism $\Phi$ of \propref{semi cross} is 
suitably compatible with regular representations, we use two lemmas.
For the first, consider 
an action $\beta$ of a discrete group $G$ on a $C^*$-algebra $A$.
For any invariant closed ideal $I$ of $A$, let $q\colon A\to A/I$ be
the quotient map, and let $\tilde\beta$ be the associated action of 
$G$ on $A/I$.
Let 
$\ind q\colon A\times_\beta G\to A/I\times_{\tilde\beta,r} G$
be the unique homomorphism such that
\[
\ind q(f)=q\circ f\forall f\in C_c(G,A).
\]
Then standard techniques from \cite[Th\'eor\`eme~4.12]{ZelPC}
yield the following:

\begin{lem}\label{ind q}
With the above assumptions and notation, 
there is a unique conditional expectation
$P_{A\times G}$ of $A\times_\beta G$ onto $A$ such that
$P_{A\times G}(f)=f(e)$
for $f\in C_c(G,A)$. The composition
$q \circ P_{A \times G}$ is a conditional expectation
of $A\times_\beta G$ onto $A/I$ such that
for $b\in A\times_\beta G$,
\[
\ind q(b)=0 \midtext{if and only if} q\circ P_{A\times G}(b^*b)=0.
\]
\end{lem}

Now let $G$ act on an $r$-discrete Hausdorff groupoid $R$, and let
$\beta$ denote the action of $G$ on $C^*(R)$ such that
\[
\beta_s(f)(x)=f(s^{-1}\cdot x)
\forall f\in C_c(R),\mbox{ }s\in G,\mbox{ }x\in R. 
\]
Also let $\lambda_R\colon
C^*(R)\to C^*_r(R)$ be the regular representation, viewed as a
quotient map, and let $P_R$ be the
conditional expectation of $C^*(R)$ onto $C_0(R^0)$ such that
\[
P_R(f)=f|_{R^0}
\forall  f\in C_c(R). 
\]
Then it follows from 
\cite[Proposition~II.4.8]{RenGA} that for $b\in C^*(R)$,
$\lambda_R(b)=0$ if and only if $P_R(b^*b)=0$.

\begin{lem}\label{ker-lem}
With the above assumptions and notation, 
the kernel of the regular representation
$\lambda_R$ is a $\beta$-invariant ideal of $C^*(R)$.
\end{lem}

\begin{proof}
It suffices to show that for $b\in C^*(R)$ and $s\in G$,
$P_R(b)=0$ {if and only if} $P_R\circ\beta_s(b)=0$.
Let $f\in C_c(R)$. Then
\begin{align*}
\norm{P_R\circ\beta_s(f)}
&=\sup_{u\in R^0}\abs{\beta_s(f)(u)}
=\sup_{u\in R^0}\abs{f(s^{-1}\cdot u)}
=\sup_{u\in R^0}\abs{f(u)}
=\norm{P_R(f)}.
\end{align*}
Hence $\norm{P_R\circ\beta_s(b)}=\norm{P_R(b)}$ for all $b\in C^*(R)$,
which proves the lemma.
\end{proof}

Note that \lemref{ker-lem} ensures that the map
$\Ind\lambda_R$ is well-defined.

\begin{prop}
\label{red semi cross}
Let $G$ be a discrete group acting on an $r$-discrete
Hausdorff groupoid $R$, let $\beta$ denote the associated action
on $C^*(R)$, and let $\Phi$ be the isomorphism of \propref{semi cross}.
Then there is an isomorphism $\Phi_r$ such that the following diagram
commutes:
\begin{equation*}
\begin{diagram}
\node{C^*(R \rtimes G)}
\arrow{e,t}{\Phi}
\arrow{s,l}{\lambda_{R \rtimes G}}
\node{C^*(R)\times_\beta G}
\arrow{s,r}{\ind \lambda_R}
\\
\node{C^*_r(R \rtimes G)}
\arrow{e,b}{\Phi_r}
\node{C^*_r(R)\times_{\beta,r} G.}
\end{diagram}
\end{equation*}
\end{prop}

\begin{proof}
We need only show that
$\ker\bigl(
\ind \lambda_R \circ \Phi \bigr)
=\ker \lambda_{R\times G}$.
Take a positive element $b$ of $C^*(R \rtimes G)$. 
By \lemref{ind q},
$\ind \lambda_R \circ \Phi(b)=0$ {if and only if}
$\lambda_R \circ P_{C^*(R)\times G} \circ \Phi(b)=0$
(because $\Phi(b)$ is positive in $C^*(R)\times_\beta G$),
so that
$\ind \lambda_R \circ \Phi(b)=0$ {if and only if}
$P_R \circ P_{C^*(R)\times G} \circ \Phi(b)=0$
(because $P_{C^*(R)\times G} \circ \Phi(b)$ is positive in
$C^*(R)$).
On the other hand, 
$\lambda_{ R \rtimes G }(b) = 0$ if and only if $P_{ R \rtimes G }(b)=0$.
Thus, it suffices to show that for all $b\in C^*(R \rtimes G)$,
\[
\norm{ P_{R \rtimes G}(b) }
= \norm{ P_R \circ P_{C^*(R)\times G} \circ \Phi(b) },
\]
and for this it suffices to take $b\in C_c(R \rtimes G)$:
\begin{gather*}
\norm{ P_R \circ P_{ C^*(R)\times G } \circ \Phi(b) }
= \norm{ P_{ C^*(R)\times G } \circ \Phi(b) |_{R^0} }
= \norm{ \Phi(b)(e) |_{R^0} }
= \sup_{ u \in R^0 } \abs{ \Phi(b)(e)(u) }
\\
= \sup_{ u \in R^0 } \abs{ b(u,e) }
= \norm{ b |_{( R^0 \times \{e\} )} }
= \norm{ b |_{( R \rtimes G )^0} }
= \norm{ P_{ R \rtimes G }(b) }.
\end{gather*}
This completes the proof.
\end{proof}

We write $Q\times_c G\rtimes G$ for the semidirect
product of $G$ acting on $Q\times_c G$, and we write the elements as
triples.

\begin{prop}
\label{equiv}
Let $c$ be a continuous homomorphism of an $r$-discrete Hausdorff
groupoid $Q$ into a discrete group $G$, and let $G$ act
on the skew product $Q\times_c G$ as in Equation~\eqref{gpdact}.  
Then the semidirect-product groupoid $Q\times_c G\rtimes G$ is
equivalent to $Q$.
\end{prop}

\begin{proof}
We will show that the space $Q\times_c G$ implements a groupoid
equivalence (in the sense of \cite[Definition~2.1]{MRW-EI})
between $Q\times_c G\rtimes G$ (acting on the left) and $Q$
(acting on the right). For the right action we need a continuous open
surjection $\sigma$ from  $Q\times_c G$ onto the unit space of $Q$. For
$(x,s)\in
Q\times_c G$ define $\sigma(x,s)=s(x)$.
Then $\sigma$ is a continuous and open surjection onto $Q^0$. Now put
\[
(Q\times_c G)*Q=
\{((x,s),y)\in (Q\times_c G)\times Q\mid
\sigma(x,s)=r(y)\},
\]
and define a map $((x,s),y)\mapsto (x,s)y$ from $(Q\times_c G)*Q$ to
$Q\times_c G$ by
\[
(x,s)y=(xy,c(y)^{-1}s).
\]
The continuity and algebraic properties of this map are easily checked,
so we have a right action of $Q$ on the space $Q\times_c G$.
For the left action we need a continuous and open surjection
$\rho$ from $Q\times_c G$ onto the unit space of $Q\times_c
G\rtimes G$. Note that
this unit space is $Q^0\times G\times \{e\}$, and the range and source
maps in $Q\times_c G\rtimes G$ are given by
\[
r(x,s,t)=(r(x),c(x)s,e)
\midtext{and}
s(x,s,t)=(s(x),s,e).
\]
For $(x,s)\in Q\times_c G$ define $\rho(x,s)=(r(x),c(x)s,e)$.
Then $\rho$ is a continuous surjection onto $Q^0\times G\times \{e\}$,
and $\rho$ is open since $r$ is and $G$ is discrete. Now put
\begin{multline*}
(Q\times_c G\rtimes G)*(Q\times_c G)
\\=
\{(x,s,t),(y,r))\in (Q\times_c G\rtimes G)\times (Q\times_c G)\mid
s(x,s,t)=\rho(y,r)\},
\end{multline*}
and define a map $((x,s,t),(y,r))\mapsto (x,s,t)(y,r)$ from $(Q\times_c
G\rtimes G)*(Q\times_c G)$ by
\[
(x,s,t)(y,r)=(xy,rt^{-1}).
\]
The continuity and algebraic properties of this map are also
easily checked,
so we have a left action of $Q\times_c G\rtimes G$ on the space
$Q\times_c G$.

Next we must show that 
both actions are free and proper, and that they commute.
If $(x,s,t)(y,r)=(y,r)$, then $xy=y$ and
$rt^{-1}=r$, so $x$ is a unit and $t=e$, hence $(x,s,t)$ is a unit; thus
the left action is free. For properness of the left action,
it is enough to show that
if $L$ is compact in $Q$ and $F$ is finite in $G$, then there is some
compact set in $(Q\times_c G\rtimes G)*(Q\times_c G)$ containing
all pairs $((x,s,t),(y,r))$ for which
\[
((x,s,t)(y,r),(y,r))\in (L\times F)\times (L\times F).
\]
But the above condition forces
$x\in LL^{-1}$,
$s\in c(L)FF^{-1}F$,
$t\in F^{-1}F$,
$y\in L$, and
$r\in F$,
so the left action is proper.
Freeness and properness of the right action is checked similarly (but
more easily), and it is straightforward to verify that the actions
commute. 

To show $Q\times_c G$ is a $(Q\times_c G\rtimes G)$-$Q$
equivalence, it remains to verify that the map
$\rho$
factors through a bijection of $(Q\times_c G)/Q$ onto $(Q\times_c
G\rtimes G)^0$, and
similarly that the map $\sigma$ factors through a
bijection of
$(Q\times_c G\rtimes G)\setminus (Q\times_c G)$ onto $Q^0$.
Since $\rho$ and $\sigma$
are surjective and the actions commute, it suffices to show that
$\rho(x,s)=\rho(y,t)$ implies $(x,s)\in (y,t)Q$, and
$\sigma(x,s)=\sigma(y,t)$ implies $(x,s)\in (Q\times_c G\rtimes G)(y,t)$.
For the first, if $\rho(x,s)=\rho(y,t)$ then
$r(x)=r(y)$ and $c(x)s=c(y)t$.
Put $z=y^{-1}x$; then $x=yz$ and
$c(z)^{-1}t=c(x)^{-1}c(y)t=s$,
so $(x,s)=(y,t)z$.
For the second, if $\sigma(x,s)=\sigma(y,t)$ then
$s(x)=s(y)$. Put
$z=xy^{-1},r=c(y)s,\text{ and }q=s^{-1}t$;
then
$x=zy$, $s=tq^{-1}$, and $c(y)tq^{-1}=c(y)s=r$,
so $(x,s)=(z,r,q)(y,t)$.
\end{proof}

\begin{prop}
\label{amenact}
Let $c$ be a continuous homomorphism of an $r$-discrete Hausdorff
groupoid $Q$ into a discrete group $G$, and suppose $Q$ is amenable.  
Then the action
$\beta$ of $G$ on $C^*(Q\times_c G)$ defined by Equation \eqref{beta}
is amenable in the sense that the regular representation 
of $C^*(Q\times_cG)\times_\beta G$ is faithful.
\end{prop}

\begin{proof}
First note that
\cite[Proposition~6.1.7]{AR-AG}, for example, implies that
the full and reduced
$C^*$-algebras of an amenable groupoid coincide.  
Since $Q$ is amenable so is the skew product $Q \times_c G$,
by \cite[Proposition~II.3.8]{RenGA}; hence
$C^*( Q \times_c G ) = C^*_r( Q \times_c G )$ and
$\ind \lambda_{ Q \times_c G }$ is just the regular representation
$\lambda_{ C^*( Q \times_c G ) \times G }$.
The semidirect-product groupoid
$Q \times_c G \rtimes G$ is also amenable, by \propref{equiv}, since
groupoid equivalence preserves amenability
(\cite[Theorem~2.2.13]{AR-AG}).
Thus, \propref{red semi cross} gives a commutative diagram
\begin{equation*}
\begin{diagram}
\node{ C^*( Q \times_c G \rtimes G ) }
\arrow{e,t}{ \Phi }
\arrow{se,b}{ \Phi_r }
\node{ C^*( Q \times_c G ) \times_\beta G }
\arrow{s,r}{ \lambda_{ C^*( Q \times_c G ) \times G } }
\\
\node[2]{ C^*( Q \times_c G ) \times_{\beta,r} G }
\end{diagram}
\end{equation*}
in which $\Phi$ and $\Phi_r$ are isomorphisms.  
This proves the proposition.
\end{proof}

\begin{rem}
The above result could also be proved using \cite[Th\'eor\`eme 4.5 and
Proposition 4.8]{AnaSD}, since both $C^*(Q \times_c G)$ and $C^*(Q
\times_c G)\times_{\beta,r} G$ are nuclear
(by \cite[Proposition~3.3.5 and Corollary~6.2.14]{AR-AG} and
\cite[Proposition II.3.8]{RenGA}).
\end{rem}

\begin{lem}\label{amen-cond-lem}
Let $c$ be a continuous homomorphism of an $r$-discrete Hausdorff
groupoid $Q$ into a discrete group $G$, and put $N=c^{-1}(e)$.  Assume
that $Q$ is amenable.  Then the canonical map $i\colon C^*(N)\to
C^*(Q)$ is faithful. 
\end{lem}

\begin{proof}
Since $Q$ is amenable, so is $N$ \cite[Proposition~5.1.1]{AR-AG}.
Let $P_Q\colon C^*(Q)\to C_0(Q^0)$ denote the unique conditional
expectation extending the map $f\mapsto f|_{Q^0}$ at the level of
$C_c$-functions. Since $Q$ is amenable, the regular representation of
$C^*(Q)$ onto $C^*_r(Q)$ is faithful \cite[Proposition~6.1.7]{AR-AG}.
By \cite[Proposition~II.4.8]{RenGA}, this implies
$P_Q$ is faithful in the sense
that $a\in C^*(Q)$ and $P_Q(a^*a)=0$ imply $a=0$, and similarly for
$P_N$ (Renault assumes $Q$ is principal, but this is not used in
showing
his conditional expectation
is faithful on the reduced $C^*$-algebra $C^*_r(Q)$).
It is easy to see by checking elements of $C_c(N)$ that $P_N=P_Q\circ
i$.  If $a\in \ker i$ then so is $a^*a$, thus $P_N(a^*a)=0$, so
$a^*a=0$ since $N$ is amenable, hence $a=0$.
\end{proof}

It is easy to check that roughly the same
argument as above would work if we only assume $N$ itself is amenable.

\begin{thm}
\label{full-gpd-thm}
Let $c$ be a continuous homomorphism of an amenable $r$-discrete Hausdorff 
groupoid $Q$ into a discrete group $G$, and let $\beta$ be the action
of Equation~\eqref{beta}. Then
\[
C^*(Q\times_c G)\times_\beta G\cong C^*(Q)\otimes\c K(\ell^2(G)).
\]
\end{thm}

\begin{proof}
\lemref{amen-cond-lem} ensures that the hypotheses of
\corref{duality} are satisfied, which gives
\[
C^*(Q\times_c G)\times_{\beta,r} G\cong C^*(Q)\otimes\c K(\ell^2(G)).
\]
The theorem now follows from \propref{amenact}.  
\end{proof}


\section{Embedding $C^*(N)$ in $C^*(Q)$}
\label{embed}

In this section we fulfill the promise made just before 
\thmref{gpdiso} by showing the map $i\colon C^*(N)\to C^*(Q)$ is
faithful when $Q$ is second countable. But first we need the following
elementary lemma, which we could not find in the literature.

\begin{lem}
Let $Q$ be an $r$-discrete Hausdorff groupoid, and let $\pi$ be a
$^*$-homomorphism from $C_c(Q)$ to the $^*$-algebra of adjointable
linear operators on an inner product space $\c H$. Then for all $a\in
C_c(Q)$, the operator $\pi(a)$ is bounded and
$\norm{\pi(a)}\le\norm{a}$, where $C_c(Q)$ is given the largest
$C^*$-norm.
\end{lem}

\begin{proof}
Let $a\in C_c(Q)$. Since $C_c(Q)$ has the largest $C^*$-norm, it
suffices to show $\pi(a)$ is bounded. Choose open bisections
(``$Q$-sets'', in Renault's terminology)
$\{U_i\}_1^n$ of $Q$ such that $\supp a\subseteq\bigcup_1^n U_i$, and a
partition of unity $\{\phi_i\}_1^n$ subordinate to the open cover
$\{U_i\}_1^n$ of $\supp a$. Then $a=\sum_1^n a\phi_i$, and $\supp
a\phi_i\subseteq U_i$. Conclusion: without loss of generality there
exists an open bisection $U$ of $Q$ such that $\supp a\subseteq U$. Then
$\supp a^*a\subseteq U^{-1}U$, a relatively compact subset of
the unit space $Q^0$. Choose an open set $V\subseteq Q^0$ such that
$\overline{U^{-1}U}\subseteq V$ and $\overline V$ is compact. Then
$a^*a\in C_0(V)$, which is a $C^*$-subalgebra of the commutative
$^*$-subalgebra $C_c(Q^0)$ of $C_c(Q)$. Since $\pi$ restricts to a
$^*$-homomorphism from $C_0(V)$ to the adjointable linear operators on
$\c H$, $\pi(a^*a)$ is bounded. Since $\pi(a)^*\pi(a)=\pi(a^*a)$,
$\pi(a)$ is bounded as well.
\end{proof}

\begin{thm}
\label{faith}
Let $c$ be a continuous homomorphism of an $r$-discrete Hausdorff groupoid
$Q$ into a discrete group $G$, and put $N=c^{-1}(e)$.  Assume that $Q$
is second countable. Then the canonical map $i\colon
C^*(N)\to C^*(Q)$ is
faithful.
\end{thm}

\begin{proof}
For notational convenience, throughout this proof we suppress the
map $i$, and {identify} $C_c(N)$ and $C^*(N)$ with their images in
$C^*(Q)$.  Our strategy is to find a $C^*$-seminorm on $C_c(Q)$ which
restricts to the greatest $C^*$-norm on $C_c(N)$.  This suffices, for
then \emph{a fortiori} the greatest $C^*$-norm on $C_c(Q)$ restricts to
the greatest $C^*$-norm on $C_c(N)$, which is what we need to prove.

To get this $C^*$-seminorm on $C_c(Q)$, we make $C_c(Q)$ into a
pre-Hilbert $C_c(N)$-module, and show that by left multiplication
$C_c(Q)$ acts by bounded adjointable operators.  We do this by showing
that the space
$Q$ implements a groupoid equivalence 
in the sense of \cite[Definition~2.1]{MRW-EI} between $N$ (acting on the
right) and a suitable groupoid $H$ (acting on the left); then the
construction of \cite{MRW-EI} shows that $C_c(Q)$ is a pre-imprimitivity
bimodule, and in particular a right pre-Hilbert $C_c(N)$-module.

We define
\[
H=\{(x,c(y))\mid x,y\in Q,s(x)=r(y)\},
\]
which is a subgroupoid of the skew product $Q\times_c G$. We claim that
$H$ is open in $Q\times_c G$.  Let $(x,t)\in H$.  There exists
$y\in Q^{s(x)}$ such that $c(y)=t$, and then there exists a
neighborhood $V$ of $y$ such that $c(V)\subseteq\{s\}$. Then $r(V)$ is
a neighborhood of $r(y)=s(x)$, so there exists a neighborhood $U$ of
$x$ such that $s(U)\subseteq r(V)$. By construction, for all $z\in U$
there exists $w\in V$ such that $r(w)=s(z)$, and then $c(w)=t$.
Therefore, the open subset $U\times\{t\}$ of $Q\times G$ is contained
in $H$, so $(x,t)$ is an interior point of $H$.  This proves
the claim.  Since the restriction of a Haar system to an open
subgroupoid gives a Haar system, counting measures give a Haar system
on $H$. Since $Q$ is second countable, the image of the
homomorphism $c$ in $G$ is countable, hence the groupoid $H$ is
second countable. Since the skew-product groupoid $Q\times_c G$ is
$r$-discrete, so is the open subgroupoid $H$.

The subgroupoid $N$ acts on the right of $Q$ by multiplication.  We
want to define a left action of the groupoid $H$ on the
space $Q$. For this we need a continuous and open surjection $\rho$
from $Q$ onto the unit space of $H$. We have
$H^0=\{(u,t)\in Q^0\times G\mid t\in c(Q^u)\}$,
and the range and source maps in $H$ are given by
\[
r(x,t)=(r(x),c(x)t)\midtext{and}s(x,t)=(s(x),t).
\]
For $y \in Q$ define
\[
\rho(y)=(r(y),c(y)).
\]
Then $\rho$ is a continuous surjection onto $H^0$, and is
open since $r$ is and $G$ is discrete.  Now put
\[
 H*Q=\{((x,t),y)\in  H\times Q\mid s(x,t)=\rho(y)\},
\]
and define a map $((x,t),y)\mapsto (x,t)y$ from $ H*Q$ to
$Q$ by
\[
(x,t)y=xy.
\]
The continuity and algebraic properties of this map are easily
checked, so we have an action of $H$ on $Q$.

Next we must show that
both actions are free and proper, and that the actions commute.
Since
$(x,t)y=y$ implies that $x$, hence $(x,t)$, is a unit, the left action
is free. For properness of the left action, let $K$ be a compact subset
of $Q\times Q$. We must show that the inverse image of $K$ under the
map $((x,t),y)\mapsto ((x,t)y,y)$ from $ H*Q$ to $Q\times Q$ is
compact. Without loss
of generality suppose $K=L\times L$ for some compact subset $L$ of $Q$.
For all $((x,t),y)\in H*Q$, if $((x,t)y,y)\in L\times L$ then
$x\in LL^{-1}$, $t\in c(L)$, and $y\in L$,
so the inverse image of $K$ is contained in
\[
(LL^{-1}\times c(L))*L,
\]
which is compact in $(Q\times_c G)*Q$.
It is easier to see that the right $N$-action is free and proper, and
straightforward to check that the actions commute. 

To show $Q$ is an $H$--$N$
equivalence, it remains to verify that the map $\rho$
factors through a bijection of $Q/N$ onto $H^0$, and
similarly that the map $s\colon Q\to H^0$ factors through a
bijection of $ H\backslash Q$ onto $N^0$.  Since $\rho$ and $s$
are surjective and the actions commute, it suffices to show that
$\rho(y)=\rho(z)$ implies $z\in yN$ and 
$s(y)=s(z)$ implies $z\in  Hy$.
For the first, if $\rho(y)=\rho(z)$ then $r(y)=y(z)$ and $c(y)=c(z)$.
Put $n=y^{-1}z$. Then $c(n)=c(y)^{-1}c(z)=e$, so $n\in N$, and $z=yn$.
For the second, if $s(y)=s(z)$, put $x=zy^{-1}$.
Then $(x,c(z))\in H$ and $z=xy=(x,c(z))y$.

Now the theory of
\cite{MRW-EI} tells us $C_c(Q)$ becomes a pre-Hilbert $C_c(N)$-module,
where $C_c(N)$ is given the $C^*$-norm from $C^*(N)$.
{From} the formulas in \cite{MRW-EI} the right module multiplication is
given by
\begin{equation*}
ac(x)=\sum_{r(n)=s(x)}a(xn)c(n^{-1}),
\end{equation*}
where $a\in C_c(Q)$ and $c\in C_c(N)$, and the inner product is
\begin{equation}
\label{inner}
\rip{a,b}{C_c(N)}(n)=\sum_{r(x,s)=\rho(y)}
\overline{a((x,s)^{-1}y)}b((x,s)^{-1}yn),
\end{equation}
where $a,b\in C_c(Q)$ and
$y$ is any element of $Q$ with $s(y)=r(n)$. The right module
action is just right multiplication by the subalgebra $C_c(N)$ inside
the algebra $C_c(Q)$. The inner product also simplifies in our
situation:
let $a,b\in C_c(Q)$, and write $a=\sum_{t\in G}a_t$ and $b=\sum_{t\in
G}b_t$ with $a_t,b_t\in C_t = \{f\in C_c(Q)\mid \supp f\subseteq
c^{-1}(t)\}$. We claim that
\[
\rip{a,b}{C_c(N)}=\sum_{t\in G}a_t^*b_t.
\]
Of course, we are identifying $a_t^*b_t$ with $a_t^*b_t|_N$, but this
causes no harm since $a_t^*b_t$ is supported in $N$. 
In Equation~\eqref{inner} we can take $y=r(n)$, so that
$\rho(y)=(r(n),e)$. Then the condition $r(x,s)=\rho(y)$ becomes
$r(x)=r(n)$ and $c(x)s=e$, so that
\begin{align*}
\rip{a,b}{C_c(N)}(n)
&=\sum_{\substack{r(x)=r(n)\\c(x)s=e}}
\overline{a((x^{-1},c(x)s)r(n))}
b((x^{-1},c(x)s)n)\\
&=\sum_{\substack{r(x)=r(n)\\c(x)s=e}}
\overline{a(x^{-1})}b(x^{-1}n)\\
&=\sum_{\substack{r(x)=r(n)\\c(x)s=e}}
a^*(x)b(x^{-1}n)\\
&=\sum_{t,r\in G}\sum_{\substack{r(x)=r(n)\\c(x)s=e}}
a^*_t(x)b_r(x^{-1}n)\\
&=\sum_{t\in G}\sum_{r(x)=r(n)}a^*_t(x)b_t(x^{-1}n).
\end{align*}
Since in this last expression
we need only consider terms with $c(x)=t^{-1}$ and
$c(x^{-1}n)=r$, which forces $t=r$, and then $s=t$ in the inner sum,
this gives
\[
\rip{a,b}{C_c(N)}(n)
=\sum_{t\in G}a^*_tb_t(n).
\]
This proves the claim.

Now we show that for fixed $a\in C_c(Q)$, the map
$b\mapsto ab$ is a 
bounded adjointable operator on the pre-Hilbert $C_c(N)$-module
$C_c(Q)$, with adjoint $b\mapsto a^*b$. This will give a
representation of $C_c(Q)$ in $\c L_{C_c(N)}(C_c(Q))$,
hence a $C^*$-seminorm on $C_c(Q)$.

We first handle the adjointability. Without loss of generality let
$a\in C_s$ and take
$b=\sum_tb_t,c=\sum_tc_t\in C_c(Q)$ with $b_t,c_t\in C_t$. Then
\begin{align*}
\rip{ab,c}{C_c(N)}
&=\rip{\textstyle\sum_t ab_t,\sum_t c_t}{C_c(N)}
=\textstyle\sum_t (ab_t)^*c_{st}
\righttext{(since $ab_t\in C_{st}$)}\\
&=\textstyle\sum_t b^*_ta^*c_{st}
=\rip{\textstyle\sum_t b_t,\sum_t a^*c_{st}}{C_c(N)}
\righttext{(since $a^*c_{st}\in C_t$)}\\
&=\rip{b,a^*\textstyle\sum_t c_{st}}{C_c(N)}
=\rip{b,a^*c}{C_c(N)}.
\end{align*}

For the boundedness, let $\omega$ be a state on $C^*(N)$, and
let $\rip{\cdot,\cdot}\omega =\omega(\rip{\cdot,\cdot}{C_c(N)})$
be the associated semi-inner product on $ C_c(Q)$. Let $\c H$
be the corresponding inner product space, and let $\Theta\colon
 C_c(Q)\to \c H$ be the quotient map. Then left multiplication
defines a $^*$-homomorphism $\pi$
from $C_c(Q)$ to the $^*$-algebra of adjointable linear
operators on $\c H$ via
$\pi(a)\Theta(b)=\Theta(ab)$.
As we show in the general lemma below,
for all $a\in C_c(Q)$,
the operator $\pi(a)$ is bounded and $\norm{\pi(a)}\le\norm{a}$.
Hence, for all $a\in C_c(Q)$ and $b\in C_c(Q)$,
\begin{multline*}
\omega\bigl(\rip{ a b,a b}{C_c(N)}\bigr)
=\rip{\pi(a)\Theta(b),\pi(a)\Theta(b)}\omega\\
\le\norm{\pi(a)}^2\rip{\Theta(b),\Theta(b)}\omega
\le\norm{a}^2\omega\bigl(\rip{b,b}{C_c(N)}\bigr).
\end{multline*}
Since the state $\omega$ was arbitrary,
\[
\rip{a b,a b}{C_c(N)}\le\norm{a}^2\rip{b,b}{C_c(N)},
\]
as required.

We can now define a $C^*$-seminorm $\norm{\cdot}_*$ on $ C_c(Q)$ by
letting $\norm{a}_*$ be the norm of the operator $b\mapsto ab$ in $\c
L_{C_c(N)}(C_c(Q))$.  To finish, we need to know that for $a\in C_c(N)$
the norm $\norm{a}_*$ agrees with the greatest $C^*$-norm
$\norm{a}$,
and it suffices to show $\norm{a}\le
\norm{a}_*$:
\[
\norm{a}^2=\norm{a^*a}\le\norm{a^*}_*\norm{a}=\norm{a}_*\norm{a},
\]
since $a^*a$ is a value of the operator $c\mapsto a^*c$, and then
canceling $\norm{a}$ gives the desired inequality.
This completes the proof. 
\end{proof}


\providecommand{\bysame}{\leavevmode\hbox to3em{\hrulefill}\thinspace}

\end{document}